\documentclass[pre-print,authoryear]{elsarticle}
\usepackage{amssymb,amsthm,amsmath,graphicx,color,textcomp}
\usepackage{amsfonts}
\usepackage{epsfig,latexsym,graphicx}
\usepackage{caption}
\usepackage{subfig}
\newtheorem{theorem}{Theorem}[section]

\newtheorem{corollary}[theorem]{Corollary}

\theoremstyle{plain}

\biboptions{comma,round}

\begin{document}

\begin{frontmatter}

\title{Testing Composite Null Hypothesis Based on $S$-Divergences}
 \author{Abhik Ghosh\fnref{label2}}
 \ead{abhianik@gmail.com}
 \author{Ayanendranath Basu\fnref{label1}}
  \ead{ayanbasu@isical.ac.in}
\address{Indian Statistical Institute}
 \fntext[label2]{This is part of the Ph.D. research work of the first author
 which is ongoing at the Indian Statistical Institute}
 \fntext[label1]{Corresponding Author}

\begin{abstract}
We present a robust test for composite null hypothesis based on the general $S$-divergence family. 
This requires a non-trivial extension of the results of Ghosh et al.~(2015).
We derive the asymptotic and theoretical robustness properties of the resulting test
along with the properties of the minimum $S$-divergence estimators under parameter restrictions 
imposed by the null hypothesis. An illustration in the context of the normal model is also presented.
\end{abstract}


\end{frontmatter}


\section{Introduction}

Statistical tests for composite hypotheses are encountered all the time in  all disciplines of applied sciences.
For such composite hypotheses, the null parameter space is generally defined through some pre-specified 
restrictions and one needs to estimate the parameter value under those restrictions to perform the test.
The most common and widely used statistical tool to solve this inferential problem is 
the classical likelihood ratio test  \citep{Neyman/Pearson:1928,Wilks:1938} 
which utilizes the maximum likelihood estimator of the parameters under given restrictions. 
However, the non-robust nature of such likelihood based solutions 
under misspecification of models and/or presence of outliers is well-known.
So, there have been many attempts for developing robust alternative to the 
likelihood ratio test (LRT) with good asymptotic and robustness properties.

\cite{Ghosh/etc:2015} proposed a general family of tests of hypothesis for the simple null problem;
these tests are based on the family of $S$-divergences \citep{Ghosh/etc:2013a},
and extend the idea of \cite{Basu/etc:2013a}
who considered testing of hypothesis based on the density power divergence (Basu et al., 1998).
In the present paper we provide a non-trivial generalization of the \cite{Ghosh/etc:2015} paper, and 
the theoretical robustness properties described in this work also provide
the theoretical underpinnings of the \cite{Basu/etc:2013b} tests as a special case. 
We also study the corresponding minimum divergence estimators under the restrictions imposed by the null.
%

In this paper, we presents the restricted minimum $S$-divergence estimator 
and its asymptotic distribution for both the discrete and continuous models.
The main focus of the paper is on the theoretical robustness properties of the 
$S$-divergence based tests of composite hypothesis. 
For brevity in presentation, the proofs of all the results are
provided in the online supplement to this paper.

%

\section{The Restricted Minimum $S$-Divergence Estimators (RMSDE)}
\label{SEC:6RMSDE}
  
The $S$-divergence family has been recently introduced by \cite{Ghosh/etc:2013a}
and contains several popular density-based divergences like the power divergence (PD) family of 
\cite{Cressie/Read:1984} and the density power divergence (DPD) family of \cite{Basu/etc:1998}. 
For two densities $g$ and $f$, it is defined in terms of two parameters $\gamma \in [0, 1]$ 
and $\lambda \in \mathbb{R}$ as
\begin{equation}
S_{(\gamma, \lambda)}(g,f) =  \frac{1}{A} ~ \int ~ f^{1+\gamma}  -   \frac{1+\gamma}{A B} ~ 
\int ~~ f^{B} g^{A}  + \frac{1}{B} ~ \int ~~ g^{1+\gamma}, ~~~~~~ A \ne 0,~B \ne 0,
\label{EQ:S_div_gen}
\end{equation}
where $A = 1+\lambda (1-\gamma)$ and  $B = \gamma - \lambda (1-\gamma)$.
Whenever $A=0$ or $B=0$, the corresponding $S$-divergence measure is defined by 
the continuous limits of (\ref{EQ:S_div_gen}) as $A\rightarrow0$ or $B\rightarrow0$ respectively.
Several properties and applications of the (unrestricted) minimum $S$-divergence estimators
have been studied by \cite{Ghosh/etc:2013a,Ghosh/etc:2013},
\cite{Ghosh:2014b}, \cite{Ghosh/Basu:2014} and \cite{Ghosh/etc:2015}.
Here, we consider the minimum $S$-divergence estimators under some pre-specified parameter restrictions 
and study their asymptotic properties.
The general theory of robustness for general minimum divergence estimators 
under parameter restrictions has been recently developed by \cite{Ghosh:2014a},
which also contains the case of the $S$-divergence measures.

\subsection{Definition and Estimating Equation}

Consider the standard set-up of parametric inference, where we have a sample $X_1, \ldots, X_n$
from the true density $g$ which is modeled by a parametric family of densities 
$\mathcal{F}=\{f_\theta: \theta \in \Theta\subseteq\mathbb{R}^p\}$.
We assume a set of $r$ restrictions on the parameter $\theta$  given by 
\begin{equation}\label{EQ:6restrictions}
  h(\theta)=0,
\end{equation}
such that the $p\times r$ matrix defined by
$H(\theta)= \frac{\partial h(\theta)}{\partial \theta}$
exists with rank $r$ and is a continuous function of $\theta$. 

The restricted minimum $S$-divergence estimator (RMSDE) of $\theta$ is to be obtained by minimizing
$S_{(\alpha, \lambda)}(\hat{g},f_\theta)$ subject to the constraints (\ref{EQ:6restrictions}); 
here $\hat{g}$ is some non-parametric estimator of the true density $g$; this is given by the 
relative frequency of the values in the sample space for discrete models and 
by some kernel density estimator for continuous models. 
See \cite{Ghosh:2014b} and \cite{Ghosh/Basu:2014} for corresponding descriptions in the unrestricted case.
Then, using the method of Lagrange multipliers, the estimating equation of the RMSDE is given by 
\begin{eqnarray}
 \left.\begin{array}{rcl}
    \int K(\delta(x))f_{\theta}^{1+\alpha}(x) u_{\theta}(x)dx  + H(\theta)\lambda_n &=& 0 \\ 
    \mbox{and} ~~~~ h(\theta) &=& 0 
 \end{array}\right\},
    \label{EQ:6RMSDE_est_equation}
\end{eqnarray}
where $\delta(x) = \delta_n(x) = \frac{\hat{g}(x)}{f_{\theta}(x)} - 1$, 
$ K(\delta) = \frac{[(\delta+1)^A - 1]}{A}$ with $A=1+\lambda(1-\alpha)$, 
$\lambda_n$ is the vector of Lagrange multipliers, and 
$u_\theta=\frac{\partial}{\partial\theta}\log f_\theta$.

\subsection{Asymptotic Distribution under Discrete Models}
\label{SEC:6asymp_RMSDE_discrete}

First we consider the case of discrete distributions, 
where both $g$ and $f_\theta$ are densities with respect to some counting measure 
over the support $\chi = \{0, 1, 2, \cdots \}$. 
We  denote the relative frequency at any point $x$
by $r_n(x)$. Then the RMSDE can be obtained as the solution of the 
estimating equation (\ref{EQ:6RMSDE_est_equation}) with $\hat{g}(x)$ replaced by $r_n(x)$ and 
the integral replaced by the countable sum over the support $\chi$.

To prove the asymptotic properties of the RMSDE under this set-up,
we define $\widetilde{\theta^g} = \displaystyle\arg\min_{\theta: h(\theta)=0} S_{(\alpha, \lambda)}(g,{f}_\theta)$, 
the restricted ``best fitting parameter" under $g$.
Then, we have the following result under the conditions (SA1)--(SA7) of Ghosh (2014b)
provided in the Supplementary material.

\begin{theorem}\label{THM:6asymp_RMSDE_discrete}
Consider the above set-up of discrete models and assume that the conditions 
(SA1)--(SA7) of Ghosh (2014b) hold with respect to $\Theta_0 = \{\theta : h(\theta)=0\}$ 
(instead of $\Theta$). Then, we have the following:
	\begin{enumerate}
	\item [(i)] There exists a consistent sequence $\widetilde{\theta}_n$ of roots to the restricted
 minimum $S$-divergence  estimating equations (\ref{EQ:6RMSDE_est_equation}).
	\item [(ii)] Asymptotically, 
	$\sqrt{n} \left(\widetilde{\theta}_n - \widetilde{\theta}^g\right)\sim 
	N_p\left(0, \widetilde{P_g} \widetilde{V_g} \widetilde{P_g}\right),$
	where the matrices $\widetilde{P_g}$ and $\widetilde{V_g}$ are as defined in Definition 1.1
	in the Supplementary material.
\end{enumerate}
\end{theorem}
%


Next consider the particular case of the model density with $g = f_{\theta_0}$ for some $\theta_0 \in \Theta$ 
satisfying the given restriction (\ref{EQ:6restrictions}). In this case, 
we have $\widetilde{\theta^g} = \theta_0$ and hence 
$
\sqrt{n} \left(\widetilde{\theta}_n - \widetilde{\theta}^g\right) \displaystyle\mathop{\rightarrow}^\mathcal{D} 
N\left(0, ~ \widetilde{P_\alpha}(\theta_0)\widetilde{V_\alpha}(\theta_0)\widetilde{P_\alpha}(\theta_0) \right)
$
asymptotically, where $\widetilde{P_\alpha}(\theta_0)$ and $\widetilde{V_\alpha}(\theta_0)$
are defined in Definition 1.1 in the Supplementary material.
Interestingly the asymptotic distribution of the RMSDE at the model 
is also independent of the parameter $\lambda$ defining the $S$-divergence measure --- just as in the 
case of the unrestricted MSDE. It also coincides with the asymptotic distribution of the 
restricted minimum DPD estimators as obtained in  \cite{Basu/etc:2013b} independently.

\subsection{The Basu--Lindsay Approach for the RMSDE under Continuous Models}
\label{SEC:6asymp_RMSDE_continuous}

Now we consider the case of continuous models, 
where the densities $g$ and $f_\theta$ are both continuous with respect to some common dominating measure.
However, 
there is a clear incompatibility of measures between the data that are discrete and 
the assumed continuous model; hence we need to use kernel density estimator in place of $\hat{g}$ 
in the estimating equation (\ref{EQ:6RMSDE_est_equation}) which  
brings in several complications like bandwidth selection, curse of dimensionality etc.~and 
complicated conditions are needed for the asymptotic results.
The approach of Basu--Lindsay \citep{Basu/Lindsay:1994} helps us to avoid such complications by 
using convolution of the assumed model also by the same kernel; 
see \cite{Ghosh/Basu:2014} for several advantages of this approach and corresponding derivations 
in the unrestricted case.

Let us define the kernel density estimator $g^*_n$ and the corresponding smoothed versions $g^*$ 
and $f^*_\theta$ of the densities $g$ and $f_\theta$ respectively: 
\begin{eqnarray}
	g_n^*(x) &=& \int ~ W(x,y,h_n) dG_n(y) = \frac{1}{n} ~ \sum_{i=1}^n ~W(x,X_i,h_n),\nonumber\\
	g^*(x)&=&\int W(x,y,h)~ dG(y),~~~~\mbox{ and }~~~
	f_\theta^*(x)=\int W(x,y,h)~ dF_\theta(y),\nonumber
\end{eqnarray}
where $W(x, y, h_n)$ is a smooth kernel function with bandwidth $h_n$, 
$G_n$ is the empirical distribution function and $G$, $F_\theta$ are distribution functions 
of $g$ and $f_\theta$ respectively. 
Using the Basu--Lindsay approach, the restricted 
minimum divergence estimator 
is to be obtained by minimizing the $S$-divergence between $g^*_n$ and $f^*_\theta$,
subject to the restriction (\ref{EQ:6restrictions}). 
Thus, the corresponding estimating equation is given by 
\begin{eqnarray}
\left.  \begin{array}{r c l}
 \int K(\delta_n^*(x)) f_\theta^*(x)^{1+\alpha} \widetilde{u_\theta}(x) dx  + H(\theta)\lambda_n &=& 0 \\
  h(\theta)  &=& 0   
\end{array} 
\right\}.
\label{EQ:6RMSDE*_est_equation}
  \end{eqnarray}
where $\delta_n^*(x) = \frac{g_n^*(x)}{f_\theta^*(x)} -1$ and 
$\widetilde{u_{\theta}}(x) = \nabla \log f_\theta^*(x)$.
In general, the resulting estimator is not the same as the RMSDE obtained 
by minimizing $S_{(\alpha,\lambda)}(g_n^*,f_\theta)$ over $\Theta_0$; 
we denote it as the restricted minimum $S^*$-divergence estimator ($\textrm{RMSDE}^*$).
%
%
We follow Ghosh and Basu (2014) to derive the asymptotic distribution of the $\textrm{RMSDE}^*$.
Let  $\widetilde{\theta^g}^*=\displaystyle\arg\min_{\theta\in\Theta_0} S_{(\alpha, \lambda)}(g^*,{f}_\theta^*)$ be
the restricted ``best fitting parameter" under $g$.
%
Then we have the following theorem.

\begin{theorem}\label{THM:6asymp_RMSDE_cont}
Consider the above set-up of continuous models and assume that the conditions 
(SB1)--(SB7) of Ghosh and Basu (2014), 
presented in the Supplementary material, hold with respect to 
$\Theta_0 = \{\theta : h(\theta)=0\}$.
Then, 
	\begin{enumerate}
	\item [(i)] There exists a consistent sequence $\widetilde{\theta}_n^*$ of roots to the restricted
 minimum $S^*$-divergence  estimating equations (\ref{EQ:6RMSDE*_est_equation}).
	\item [(ii)] Asymptotically,
	$\sqrt{n} \left(\widetilde{\theta}_n^* - \widetilde{\theta^g}^*\right) \sim 
	N_p\left(0,\widetilde{P_{\alpha, \lambda}}^*(\widetilde{\theta^g}^*) \widetilde{V_{\alpha, \lambda}}^*(g)  \widetilde{P_{\alpha, \lambda}}^*(\widetilde{\theta^g}^*)\right)$, where 
	$\widetilde{P_{\alpha, \lambda}}^*(\widetilde{\theta^g}^*)$ and  $\widetilde{V_{\alpha, \lambda}}^*(g)$
	are defined in Definition 1.2 in the Supplementary material.
\end{enumerate}
\end{theorem}
%


Now, for $g = f_{\theta_0}$ with $\theta_0 \in \Theta$ 
satisfying (\ref{EQ:6restrictions}), 
we have $\widetilde{\theta^g}^* = \theta_0$ and 
the asymptotic distribution of 
$\sqrt{n} \left(\widetilde{\theta}_n^* - \theta_0 \right)$ 
is normal with mean 0 and variance 
$\widetilde{P_\alpha}^*(\theta_0)\widetilde{V_\alpha}^*(\theta_0)\widetilde{P_\alpha}^*(\theta_0)$,
 where $\widetilde{P_\alpha}^*(\theta_0)$ and  $\widetilde{V_\alpha}^*(\theta_0)$ are as in Definition 1.2 in the 
 Supplementary material. 
Once again, the asymptotic distribution of the RMSDE$^*$
turns out to be independent of the parameter $\lambda$ at the model.

Further,  if we assume that the kernel used in smoothing is $\alpha$-transparent
for the restricted family ${\cal F}_0 = \{f_\theta: \theta \in \Theta_0 \} \subset \mathcal{F}$ 
in the sense of Definition 9.1 of Ghosh and Basu (2014) presented in the Supplementary material (Definition 1.3),
then it follows that  for any $\theta_0 \in \Theta_0$, 
$\widetilde{J_\alpha}^*(\theta_0)=\widetilde{J_\alpha}(\theta_0)$, 
$\widetilde{V_\alpha}^*(\theta_0)=\widetilde{V_\alpha}(\theta_0)$ and so 
$\widetilde{P_\alpha}^*(\theta_0)=\widetilde{P_\alpha}(\theta_0)$. 
Hence at the model density $g=f_{\theta_0} \in \mathcal{F}_0$, the  asymptotic 
distribution of the $\textrm{RMSDE}^*$ becomes exactly the same as that of the RMSDE under discrete model. 

\section{$S$-Divergence based Test (SDT) for Composite Hypothesis}
\label{SEC:7sample1_compositeTest}

Now we consider the problem of testing composite null hypothesis.  
Under the notations of the previous section, 
take a fixed (proper) subspace $\Theta_0$ of the parameter space $\Theta$.  
Our objective is to test for the hypothesis 
\begin{equation}\label{EQ:7composite_hypo}
H_0 : \theta \in \Theta_0 ~~~ \mbox{against} ~~~~ H_1 : \theta  \notin \Theta_0.
\end{equation}

Following  Ghosh et al.~(2015), the $S$-divergence based test (SDT) statistics  
for testing the above hypothesis can be constructed as
\begin{eqnarray}
\widetilde{T_{\gamma,\lambda}^{(1)}}({\widehat{\theta}_{\beta,\tau}}, \widetilde{\theta}_{\beta,\tau}) 
= 2 n S_{(\gamma,\lambda)}(f_{\widehat{\theta}_{\beta,\tau}}, f_{\widetilde{\theta}_{\beta,\tau}}),
\label{EQ:7one_composite_TS}
\end{eqnarray}
where $S_{(\gamma,\lambda)}(\cdot, \cdot)$ is the $S$-divergence measure with parameter $\gamma$ and $\lambda$
and $\widetilde{\theta}_{\beta,\tau}$, $\widehat{\theta}_{\beta,\tau}$ denote the 
restricted (under $\Theta_0$) and unrestricted MSDEs with tuning parameters $\beta$ and $\tau$.
Their asymptotic distributions at the model are independent of $\tau$.
Indeed the estimators $\widehat{\theta}_{\beta,\tau_1}$ and $\widehat{\theta}_{\beta,\tau_2}$
are asymptotically equivalent in that 
$$
\sqrt{n}\left(\widehat{\theta}_{\beta,\tau_1} - \widehat{\theta}_{\beta,\tau_2}\right)
\mathop{\rightarrow}^\mathcal{P} 0,
$$
for any $\tau_1 \neq \tau_2$. We therefore replace $\widehat{\theta}_{\beta,\tau}$
and $\widetilde{\theta}_{\beta,\tau}$ by $\widehat{\theta}_{\beta,0}$ and $\widetilde{\theta}_{\beta,0}$
respectively, without altering the asymptotic properties of the test statistics in (\ref{EQ:7one_composite_TS}).
The advantage of this substitution is that the latter set of estimators minimizes the DPD, 
and thus can be evaluated without any kernel smoothing.
The asymptotic properties of the restricted MDPDE 
$\widetilde{\theta}_{\beta} =\widetilde{\theta}_{\beta,0}$ at the model $g=f_{\theta_0}$
is given by, see \cite{Basu/etc:2013b},
\begin{eqnarray}
\sqrt{n}(\widetilde{{\theta}_\beta} - \theta_0) \mathop{\rightarrow}^\mathcal{D} 
N(0, \widetilde{P_\beta}(\theta_0)\widetilde{V_\beta}(\theta_0)\widetilde{P_\beta}(\theta_0)).\nonumber
\end{eqnarray}
This distributional convergence holds under Conditions (D1)--(D5) of \cite{Basu/etc:2011} with respect to $\Theta_0$;
we refer to these 5 conditions as ``Basu et al.~conditions" throughout the rest of the paper.
We also assume the standard conditions of asymptotic inference, 
given by Assumptions A, B, C and D of \citet[][p. 429]{Lehmann:1983}; 
we refer to them as the ``Lehmann conditions". 
Both set of conditions are presented in the supplementary material.  	


Now, to explore the asymptotic properties of the proposed (modified) SDT statistics
$\widetilde{T_{\gamma,\lambda}^{(1)}}({\widehat{\theta}_{\beta}}, \widetilde{\theta}_{\beta}) 
=\widetilde{T_{\gamma,\lambda}^{(1)}}({\widehat{\theta}_{\beta,0}}, \widetilde{\theta}_{\beta,0})$
for testing the composite hypothesis (\ref{EQ:7composite_hypo}), 
we re-define the null parameter space $\Theta_0$ in terms of $r$ restrictions
of the form (\ref{EQ:6restrictions}). We also assume that the corresponding $p \times r$ matrix 
$H(\theta) = \frac{\partial h(\theta)}{\partial \theta} $
exists and it is a continuous function of $\theta$ with rank $r$. 
Indeed, this condition can be seen to hold for most parametric hypothesis.
We start with the asymptotic null distribution of the proposed SDT.

\begin{theorem}
Suppose the model density satisfies the Lehmann and Basu et al.~conditions with respect to both 
$\Theta$ and $\Theta_0$ and $H_0$ is true with $\theta_0\in\Theta_0$ being the true parameter value. 
Then, the asymptotic null distribution of the SDT statistic 	
$\widetilde{T_{\gamma,\lambda}^{(1)}}({\widehat{\theta}_\beta}, \widetilde{{\theta}_\beta})$ 
coincides with the distribution of $\sum_{i=1}^r ~  \widetilde{\zeta_i}^{\gamma, \beta}(\theta_0)Z_i^2,$
where $Z_1, \cdots,Z_r$ are independent standard normal variables, 
$\widetilde{\zeta_1}^{\gamma, \beta}(\theta_0)$, $\ldots$, $\widetilde{\zeta_r}^{\gamma, \beta}(\theta_0)$ 
are the nonzero eigenvalues of $A_\gamma(\theta_0)\widetilde{\Sigma}_\beta(\theta_0)$ with 
$$\widetilde{\Sigma}_\beta(\theta_0) 
= [J_\beta^{-1}(\theta_0) - P_\beta(\theta_0)]V_\beta(\theta_0)[J_\beta^{-1}(\theta_0) - P_\beta(\theta_0)]$$ 
and $r = rank\left(V_\beta(\theta_0)[J_\beta^{-1}(\theta_0) - P_\beta(\theta_0)]A_\gamma(\theta_0)
[J_\beta^{-1}(\theta_0) - P_\beta(\theta_0)] V_\beta(\theta_0)\right).$\\
\label{THM:7asymp_null_composite}
\end{theorem}

Noting the similarity of the above asymptotic null distribution with the case of 
testing simple null hypothesis,  we can find critical values of the proposed SDT
following Remark 3 of \cite{Basu/etc:2013a}.
We can also derive an asymptotic power approximation 
at any point $\theta^* \notin \Theta_0$; if $\theta^* \notin \Theta_0$ is the true parameter value 
then $\widehat{\theta}_\beta \displaystyle\mathop{\rightarrow}^\mathcal{P} \theta^*$ but
$\widetilde{{\theta}_\beta} \displaystyle\mathop{\rightarrow}^\mathcal{P} \theta_0$ for some 
$\theta_0 \in \Theta_0$ with $\theta^* \neq \theta_0$. 
Define $\Sigma_\beta(\theta)=J_\beta^{-1}(\theta)V_\beta(\theta)J_\beta^{-1}(\theta) $. 
Then, by an argument similar to the one in the above theorem, one can show under the 
{\it Basu et al.~conditions} that 
\begin{eqnarray}
  \sqrt{n}  \begin{pmatrix}
      ~~\widehat{\theta}_\beta  - \theta^* ~\\
      ~~\widetilde{{\theta}_\beta}  - \theta_0~
     \end{pmatrix} \mathop{\rightarrow}^\mathcal{D} N\left( \begin{bmatrix}
              ~~0~ \\
              ~~0~
             \end{bmatrix},  \begin{bmatrix}
                ~\Sigma_\beta(\theta^*) & A_{12}~\\ 
                ~A_{12}^T & P_\beta(\theta_0)V_\beta(\theta_0)P_\beta(\theta_0)~
                \end{bmatrix}\right),\nonumber
                \label{EQ:7asymp_null_nonnull}
  \end{eqnarray}
for some $p\times p$ matrix $A_{12}=A_{12}(\theta^*,\theta_0)$.
Further define 
\begin{eqnarray}
M_{1,\gamma,\lambda}(\theta^*,\theta_0) 
= \nabla S_{(\gamma,\lambda)}(f_\theta, f_{\theta_0})\big|_{\theta=\theta^*}, ~~~~~
M_{2,\gamma,\lambda}(\theta^*,\theta_0) 
= \nabla S_{(\gamma,\lambda)}(f_{\theta^*},f_{\theta})\big|_{\theta=\theta_0}.\nonumber
\end{eqnarray}

\begin{theorem}\label{THM:7asymp_power_composite}
Suppose the model density satisfies the Lehmann and Basu et al.~conditions with respect to both 
$\Theta$ and $\Theta_0$ and take any $\theta^*\notin \Theta_0$.
 An asymptotic approximation to the power function of the SDT statistic 
$\widetilde{T_{\gamma,\lambda}^{(1)}}({\widehat{\theta}_\beta}, \widetilde{{\theta}_\beta})$ 
for testing (\ref{EQ:7composite_hypo}) at the significance level $\alpha$ is given by 
\begin{eqnarray}
\widetilde{\pi_{n,\alpha}^{\beta,\gamma,\lambda}} (\theta^*) 
= 1 - \Phi \left( \frac{\sqrt{n}}{\widetilde{\sigma_{\beta,\gamma, \lambda}}(\theta^*, \theta_0)} 
\left(\frac{\widetilde{t_\alpha^{\beta,\gamma}}}{2 n} 
- S_{(\gamma,\lambda)}(f_{\theta^*}, f_{\theta_0})\right)\right), 
~~~ \theta^* \neq \theta_0,\nonumber
\end{eqnarray}
where $\widetilde{t_\alpha^{\beta,\gamma}}$ is the $(1-\alpha)^{\rm th}$ quantile of the asymptotic null
distribution of the SDT $\widetilde{T_{\gamma,\lambda}^{(1)}}({\widehat{\theta}_\beta}, \widetilde{{\theta}_\beta})$ 
and 
$~~\widetilde{\sigma_{\beta,\gamma, \lambda}}(\theta^*, \theta_0)^2 
= M_{1,\gamma,\lambda}^T \Sigma_\beta M_{1,\gamma,\lambda} 
+ M_{1,\gamma,\lambda}^T A_{12} M_{2,\gamma,\lambda} + M_{2,\gamma,\lambda}^T A_{12}^T M_{1,\gamma,\lambda} 
+ M_{2,\gamma,\lambda}^T P_\beta V_\beta P_\beta M_{2,\gamma,\lambda}.
$
\end{theorem}

The above theorem may be proved by a routine application of Taylor series and is omitted. 
This  power approximation can help us to obtain the required sample size 
in any planned experiment to achieve a desired power. 
The theorem also shows that the proposed $S$-divergence based test 
is consistent for the composite hypotheses at any $\theta^* \notin \Theta_0$.

\section{Robustness of the SDT for Composite Hypothesis}
\label{SEC:7robust_composite}


\subsection{Influence Function of the Test}
\label{SEC:7IF_test_composite}

Let us define the statistical functional corresponding to the proposed SDT
for the composite hypothesis as 
$$\widetilde{T_{\gamma,\lambda}^{(1)}}(G) 
= S_{(\gamma,\lambda)}(f_{U_{\beta}(G)},f_{\widetilde{U_{\beta}}(G)}),$$ 
where $U_{\beta}(G)$ is the MDPDE functional and $\widetilde{U_{\beta}}(G)$ 
is the restricted MDPDE functional under $\Theta_0$ as defined in \cite{Ghosh:2014a}. 
Consider the contaminated distribution  $H_\epsilon = (1-\epsilon) G + \epsilon \wedge_y$,
where $\wedge_y$ is the degenerate distribution at the contamination point $y$ 
and $\epsilon$ is the contamination proportion.
Then Hampel's first-order influence function
\citep{Hampel/etc:1986,Rousseeuw/Ronchetti:1979,Rousseeuw/Ronchetti:1981} 
of the SDT functional $\widetilde{T_{\gamma,\lambda}^{(1)}}(G) $
is given by 
\begin{eqnarray}
IF(y; \widetilde{T_{\gamma,\lambda}^{(1)}}, G) 
= \left.\frac{\partial}{\partial\epsilon}\widetilde{T_{\gamma,\lambda}^{(1)}}(H_\epsilon) 
\right|_{\epsilon=0} 
&=& M_{1,\gamma,\lambda}(U_{\beta}(G), \widetilde{U_{\beta}}(G))^T IF(y; U_{\beta}, G) \nonumber\\
&+& M_{2,\gamma,\lambda}(U_{\beta}(G), \widetilde{U_{\beta}}(G))^T IF(y; \widetilde{U_{\beta}}, G),\nonumber
\end{eqnarray}
where $IF(y; U_{\beta}, G)$ and $IF(y; \widetilde{U_{\beta}}, G)$ are the influence functions ($IF$s) of 
$U_{\beta}$ and $\widetilde{U_{\beta}}(G)$ respectively. 
Now, under the null hypothesis if $\theta_0\in \Theta_0$ is the true value of parameter 
with  $G=F_{\theta_0}$ then $U_{\beta}(F_{\theta_0}) = \theta_0$, 
$\widetilde{U_{\beta}}(F_{\theta_0}) = \theta_0$ and $M_{i,\gamma,\lambda}(\theta_0,\theta_0)=0$ 
for all $i=1, 2$; hence the first-order $IF$ of our SDT statistic 
for the composite hypothesis also becomes zero at the null.

Therefore, to assess the robustness of the test, 
we consider the second order influence function of our statistic defined as 
$IF_2(y; \widetilde{T_{\gamma,\lambda}^{(1)}}, G) = 
\left.\frac{\partial^2}{\partial^2\epsilon}\widetilde{T_{\gamma,\lambda}^{(1)}}(H_\epsilon) 
\right|_{\epsilon=0}$.
In the particular case $G=F_{\theta_0}$ with $\theta_0 \in \Theta_0$, 
this second order influence function of the SDT simplifies to  
\begin{eqnarray}
&& IF_2(y; \widetilde{T_{\gamma,\lambda}^{(1)}}, F_{\theta_0}) 
=D_\beta(y;\theta_0)^T A_\gamma(\theta_0) D_\beta(y;\theta_0),
\nonumber
\end{eqnarray}
where $D_\beta(y;\theta_0)=\left[IF(y; U_{\beta}, F_{\theta_0}) - IF(y; \widetilde{U_{\beta}}, F_{\theta_0})\right]$.
Note that the $IF$ of the SDT at the composite null is also independent of $\lambda$ 
and it is bounded if and only if  the influence function of 
the corresponding unrestricted and restricted MDPD functionals are both bounded or 
both diverge at the same rate. However, for most parametric models, these IFs 
of MDPDE and RMDPDE are seen to be bounded whenever $\beta>0$ 
but unbounded at $\beta=0$.

\subsection{Level and Power Influence Functions}
\label{SEC:7IF_power_composite}

Now we consider the influence function of level and power of the SDT for composite hypothesis. 
Since the SDT is consistent, its asymptotic power is one at any fixed alternative;
so we consider the asymptotic power under contiguous  alternatives 
$H_{1,n} : \theta_n = \theta_0 + \frac{\Delta}{\sqrt{n}} \in \Theta - \Theta_0$ 
with $\Delta \in \mathbb{R}^p - \{0\}$ and $\theta_0 \in \Theta_0$. 
Clearly, to ensure the existence of such a $\theta_0$ in $\Theta_0$ 
there must exist a limit point $\theta_0$ of the null parameter space $\Theta_0$;  
we assume $\Theta_0$ to be a closed subset of $\Theta$. 
Next, following \cite{Hampel/etc:1986}, we also consider the contaminations over 
these contiguous alternatives such that their effect tends to zero as $\theta_n$ tends to $\theta_0$ 
at the same rate to avoid confusion between the null and alternative neighborhoods 
\cite[also see][]{Huber/Carol:1970, Heritier/Ronchetti:1994, Toma/Broniatowski:2010}. 
So, consider the contaminated distributions $F_{n,\epsilon,y}^L$ and $F_{n,\epsilon,y}^L$, 
defined in Definition 1.4 of the supplementary material, for level and power respectively
and the level influence function ($LIF$) and the power influence function ($PIF$) 
as defined therein; 
also see Ghosh et al.~(2015).

Let us first derive a general expression for asymptotic power 
$\widetilde{P}(\Delta, \epsilon) = \displaystyle\lim_{n \rightarrow \infty} ~ P_{F_{n,\epsilon,y}^P}
\left(\widetilde{T_{\gamma,\lambda}^{(1)}}({\widehat{\theta}_\beta}, \widetilde{{\theta}_\beta}) 
>  \widetilde{t_\alpha^{\beta,\gamma}}\right) 
$ for testing composite hypothesis under contamination in the following theorem.
Here, $\chi_p^2$ denote a central chi-square random variable with $p$ degrees of freedom
and $\chi_{p,\delta}^2$  denote a non-central chi-square random variable with degrees of freedom $p$
and non-centrality parameter $\delta$.

\begin{theorem}
Assume that the Lehmann and Basu et al.~conditions hold for the model density and the null parameter space 
$\Theta_0$ is such that there exists a limit point $\theta_0 \in \Theta_0$ satisfying 
$\theta_n = \theta_0 + \frac{\Delta}{\sqrt{n}} \in \Theta - \Theta_0$ for all $\Delta \in \mathbb{R}^p - \{0\}$.
Then for any $\Delta \in \mathbb{R}^p$ and $\epsilon \geq 0$, we have the following:
\begin{itemize}
\item[(i)] The asymptotic distribution of 
$\widetilde{T_{\gamma,\lambda}^{(1)}}({\widehat{\theta}_\beta}, \widetilde{{\theta}_\beta})$ 
under $F_{n,\epsilon,y}^P$ is the same as that of the quadratic form 
$W^T A_\gamma(\theta_0)W$, where $W\sim N_p\left(\widetilde{\Delta^*}, \widetilde{\Sigma}_\beta(\theta_0)\right)$, where
$\widetilde{\Delta^*} = \left[ \Delta + \epsilon \left\{IF(y;U_\beta,F_{\theta_0}) 
- IF(y;\widetilde{U_\beta},F_{\theta_0})\right\}\right]$.
\\
Equivalently, this distribution is the same as that of 
$\sum_{i=1}^r ~  \widetilde{\zeta_i}^{\gamma, \beta}(\theta_0)\chi_{1,\widetilde{\delta}_i}^2,$
where
$\left(\sqrt{\widetilde{\delta}_1}, \ldots, \sqrt{\widetilde{\delta}_p}\right)^T 
=  \widetilde{V_{\beta,\gamma}}(\theta_0)\widetilde{\Sigma}_\beta^{-1/2}(\theta_0)\widetilde{\Delta^*}$
with $\widetilde{V_{\beta,\gamma}}(\theta_0)$ being the matrix of normalized eigenvectors of 
$A_\gamma(\theta_0)\widetilde{\Sigma}_\beta(\theta_0)$.

\item[(ii)] $\widetilde{P}(\Delta, \epsilon)=
\sum\limits_{v=0}^{\infty} ~ \widetilde{C_v^{\gamma, \beta}}(\theta_0, \widetilde{\Delta^*}) 
P\left(\chi_{r+2v}^2 > \frac{\widetilde{t_\alpha^{\beta,\gamma}}}{
\widetilde{\zeta_{(1)}}^{\gamma, \beta}(\theta_0)}\right)$,
where 
$\widetilde{\zeta_{(1)}}^{\gamma, \beta}(\theta_0)$ is the minimum of 
$\widetilde{\zeta_{i}}^{\gamma, \beta}(\theta_0)$s over $i=1, \ldots,r$  and 
$\widetilde{C_v^{\gamma, \beta}}(\theta_0, \widetilde{\Delta})$ is as Definition 1.5
in the Supplementary material.
\end{itemize}
%
\label{THM:7asymp_Contaminated_power_composite}
\end{theorem}

\begin{corollary}\label{COR:7contiguous_power_composite}
Putting $\epsilon=0$ in above theorem, we get the asymptotic power under the contiguous alternatives 
$H_{1,n}: \theta= \theta_n = \theta_0 + \frac{\Delta}{\sqrt{n}}$ as 
\begin{eqnarray}
\widetilde{P_0} &=& \widetilde{P}(\Delta, \epsilon=0) 
=\sum\limits_{v=0}^{\infty} ~ \widetilde{C_v^{\gamma, \beta}}(\theta_0, {\Delta}) 
P\left(\chi_{r+2v}^2 > \widetilde{t_\alpha^{\beta,\gamma}}/\widetilde{\zeta_{(1)}}^{\gamma, \beta}(\theta_0)
\right).\nonumber
\end{eqnarray}
\end{corollary}


\begin{corollary}
Putting $\Delta=0$ in above theorem, we get the asymptotic level under the probability distribution 
$F_{n,\epsilon,{y}}^L$ as 
\begin{eqnarray}
\widetilde{\alpha_\epsilon} &=&  \widetilde{P}(\Delta=0, \epsilon)  
=\sum\limits_{v=0}^{\infty} ~ \widetilde{C_v^{\gamma, \beta}}(\theta_0, \epsilon D_\beta(y, \theta_0)) 
P\left(\chi_{r+2v}^2 > \widetilde{t_\alpha^{\beta,\gamma}}/\widetilde{\zeta_{(1)}}^{\gamma, \beta}(\theta_0)\right).\nonumber
\end{eqnarray}
%
Further, if we also take $\epsilon=0$, then $F_{n,\epsilon,{y}}^L$ 
coincides with the null distribution and the asymptotic distribution of the proposed SDT obtained 
from part (i) of the above theorem coincides with asymptotic null distribution obtained independently in 
Theorem \ref{THM:7asymp_null_composite}; hence $\widetilde{\alpha_0 }= \alpha$, as expected.
\label{COR:7asymp_level_composite}
\end{corollary}

In practice, we can use finite truncation to approximate the infinite series in the 
above theorem, as discussed in Remark 3.1  of Ghosh et al.~(2015).

Finally we will compute the level and power influence function of the proposed $S$-divergence based 
test statistics for composite hypothesis from the expression of $\widetilde{P}(\Delta, \epsilon)$ 
as obtained in Theorem  \ref{THM:7asymp_power_composite}. 
In particular, the power influence function ($PIF$) comes from a simple differentiation of 
$\widetilde{P}(\Delta, \epsilon)$ at $\epsilon=0$ and then the level influence function ($LIF$) can be derived 
just by substituting $\Delta=0$. The following theorem presents the form of the $PIF$ and $LIF$ of the proposed SDT. 
Clearly, both the $LIF$ and $PIF$ can be seen to be bounded whenever the influence function of 
the MDPDE under the null and overall parameter space both are bounded or both diverges at the same rate; 
this in turn implies the size and power robustness of the proposed SDT for $\beta>0$.

\begin{theorem}\label{THM:7IF_power_composite}
Assume that the Lehmann and Basu et al.~conditions hold for the model density and 
the influence function $ IF(y;U_\beta,F_{\theta_0})$ of the minimum DPD estimator is bounded. 
Then the power and level influence function of the proposed test statistics for composite hypothesis 
have the forms
\begin{eqnarray}
PIF(y; T_{\gamma,\lambda}^{(1)}, F_{\theta_0}) &=& 
 \frac{\partial}{\partial\epsilon} \widetilde{P}(\Delta, \epsilon)\big|_{\epsilon=0}
 = D_\beta(y,\theta_0)^T C^*(\Delta,\theta_0,\gamma,\beta,\alpha)
\nonumber\\
\mbox{and  }~~~~
LIF(y; T_{\gamma,\lambda}^{(1)}, F_{\theta_0})
&=& D_\beta(y,\theta_0)^T C^*(0,\theta_0,\gamma,\beta,\alpha),\nonumber
\end{eqnarray} 
where
$C^*(\Delta,\theta_0,\gamma,\beta,\alpha)= \sum\limits_{v=0}^{\infty} ~ 
 \left[\frac{\partial}{\partial t}\widetilde{C_v^{\gamma, \beta}}(\theta_0, t)\big|_{t=\Delta}\right]
  P\left(\chi_{r+2v}^2 > \frac{\widetilde{t_\alpha^{\beta,\gamma}}}{
  \widetilde{\zeta_{(1)}}^{\gamma, \beta}(\theta_0)}\right).$
\end{theorem}

\section{Example: Testing Normal Mean with unknown variance}
\label{SEC:7SDT_normMu_unknown}


Now, we illustrate the proposed theory in the context of testing the mean of a normal distribution
 with unknown variance; the same problem with known variance is illustrated in Ghosh et al.~(2015). 
Suppose we have a sample $X_1, \ldots, X_n$ of size $n$ from a population 
having a univariate normal density $N(\mu, \sigma^2)$ with both the parameters unknown. 
Based on this sample, we want to test the hypothesis $H_0 : \mu = \mu_0$ 
for a pre-specified real number $\mu_0$; $\sigma$ is not assumed to be known.
Suppose $\widehat{\theta_\beta} = (\widehat{\mu_\beta},~\widehat{\sigma_\beta})$ is the MDPDE of 
$\theta=(\mu, ~ \sigma)$ with tuning parameter $\beta$ and the corresponding restricted MDPDE
under the above null is  $\widetilde{\theta_\beta} = ({\mu_0},~\widetilde{\sigma_\beta})$. 
Then a simple calculation shows that the proposed SDT for testing above $H_0$ is given by  
\begin{eqnarray}
\widetilde{T_{\gamma,\lambda}^{(1)}}({\widehat{\theta}_\beta}, \widetilde{\theta_\beta})  
&=&  \frac{2 n \widetilde{\kappa_\gamma}}{AB} 
\left[\frac{A}{\widehat{\sigma_\beta}^\gamma} + \frac{B}{\widetilde{\sigma_\beta}^\gamma} \right.
- \left. \frac{(1+\gamma)^{\frac{3}{2}}\widetilde{\sigma_\beta}^{1-B}\widehat{\sigma_\beta}^{1-A}}{ 
 \sqrt{B\widehat{\sigma_\beta}^2 + A\widetilde{\sigma_\beta}^2}}  
e^{-\frac{AB(\widehat{\mu}_\beta -\mu_0)^2}{2(B\widehat{\sigma_\beta}^2 
+ A\widetilde{\sigma_\beta}^2)}}\right], 
\nonumber
\end{eqnarray}
for $A, B \neq 0,$ where $\widetilde{\kappa_\gamma}=(2\pi)^{-\frac{\gamma}{2}} (1+\gamma)^{-\frac{1}{2}}$.
However, at $A=0$ or $B=0$, the above test statistic is defined in a limiting sense as 
\begin{eqnarray}
\widetilde{T_{\gamma,\lambda}^{(1)}}({\widehat{\theta}_\beta}, \widetilde{\theta_\beta})  
&=& \frac{n \widetilde{\kappa_\gamma}}{\widetilde{\sigma_\beta}^\gamma} 
\left[ \log\left(\frac{\widehat{\sigma_\beta}^2}{\widetilde{\sigma_\beta}^2}\right) 
+ \frac{1}{1+\gamma}\left(\frac{\widetilde{\sigma_\beta}^2}{\widehat{\sigma_\beta}^2} - 1\right)
+ \frac{(\widehat{\mu}_\beta -\mu_0)^2}{\widehat{\sigma_\beta}^2}\right], 
~~ A = 0,\nonumber\\
&=& \frac{n \widetilde{\kappa_\gamma}}{\widehat{\sigma_\beta}^\gamma} 
\left[ \log\left(\frac{\widetilde{\sigma_\beta}^2}{\widehat{\sigma_\beta}^2}\right) 
+ \frac{1}{1+\gamma}\left(\frac{\widehat{\sigma_\beta}^2}{\widetilde{\sigma_\beta}^2} - 1\right)
+ \frac{(\widehat{\mu}_\beta -\mu_0)^2}{\widetilde{\sigma_\beta}^2}\right], 
~~ B = 0.\nonumber
\end{eqnarray}
In particular, $\gamma=\lambda=\beta=0$, $\widetilde{\kappa_0}=1$,  
$\widehat{\theta_0}=(\bar{X},~\frac{1}{n}\sum_{i=1}^n(X_i - \bar{X})^2)$, the unrestricted MLE of $\theta$
and $\widetilde{\theta_0}=(\mu_0,~\frac{1}{n}\sum_{i=1}^n(X_i - \mu_0)^2)$, 
the restricted MLE of $\theta$ under the null hypothesis. 
Then, the SDT statistic further simplifies to 
$$\widetilde{T_{0,0}^{(1)}}({\widehat{\theta}_0}, {\theta_0})  
= n\log\left(\frac{\sum_{i=1}^n(X_i - \mu_0)^2}{\sum_{i=1}^n(X_i - \bar{X})^2}\right);
$$
this is again the likelihood ratio test statistic for the problem under consideration.
Thus, the proposed SDT is a robust generalization of the LRT.

From Theorem \ref{THM:7asymp_null_composite}, the asymptotic null distribution of the SDT statistics 
$\widetilde{T_{\gamma,\lambda}^{(1)}}({\widehat{\theta}_\beta}, \widetilde{\theta_\beta})$ 
for the composite null is also given by the distribution of $\zeta_1^{\gamma,\beta} Z_1$, 
where $Z_1\sim\chi_1^2$ and 
$\zeta_1^{\gamma,\beta}= \frac{\kappa_\gamma\upsilon_\beta}{\sigma^{\gamma+2}}$ 
with $\upsilon_\beta = \frac{(1+\beta)^3}{(1+2\beta)^{3/2}}\sigma^2$ 
being the asymptotic variance of the MDPDE $\widehat{\mu}_\beta$.
Note that, this asymptotic null distribution is exactly the same as that in case of testing 
simple hypothesis with known $\sigma$. In fact, all the asymptotic properties and the influence function 
of the test statistics for this case of composite hypothesis turns out to be exactly the same as obtained 
in the case of known $\sigma$ (Ghosh et al., 2015); 
the main reason behind this is the asymptotic independence of the MDPDE or RMDPDE of $\mu$
and $\sigma$ under the normal model. 



\subsection{A Real data Example: Telephone-fault data}

Consider an interesting real dataset containing the records on telephone line faults presented and analyzed by 
\cite{Welch:1987}; also studied by \cite{Simpson:1989b} and \cite{Basu/etc:2013a,Basu/etc:2013b}.
Table \ref{TAB:7data_telephone_fault} presents the data on the ordered differences between 
the inverse rates of test and control in 14 matched pairs of areas. 
These data could have been modeled by the normal distribution with mean $\mu$ and standard deviation $\sigma$
if not for the first observation which produces a huge outlier with respect to the remaining 13 observations.
The presence of this outlying observation changes the MLE of the parameters $\mu$ and $\sigma$ drastically
whereas the MDPDE with a slightly larger tuning parameter $\beta$ produces robust estimators 
\citep{Basu/etc:2013a,Basu/etc:2013b}; Table \ref{TAB:7MDPDE_telephone_fault} presents the MDPDE and MLE 
(it is the MDPDE with $\beta=0$) of the parameters under the full data and also after removing 
the outlying first observation.

\begin{table}[h]
\centering
\caption{Telephone-fault data}
\resizebox{\textwidth}{!}{
\begin{tabular}{l rrrrrrrrrrrrrr} \hline
Pair     &  1  & 2  & 3  & 4  & 5  & 6 & 7 & 8  & 9  & 10  & 11  & 12 & 13 & 14 \\ \hline\hline
Difference & $-988$ & $-135$ & $-78$ & 3 & 59 & 83 & 93 & 110 & 189 & 197 & 204 & 229 & 289 & 310  \\
\hline
\end{tabular}}
\label{TAB:7data_telephone_fault}
\end{table}
\begin{table}[h]
\centering
\caption{MDPDEs of $\mu$ and $\sigma$ for the Telephone-fault data}
\begin{tabular}{ll rrrrr} \hline
	&	$\beta$	&	0	&	0.05	&	0.1	&	0.2	&	0.5	\\\hline\hline
Full Data	&	$\widehat\mu$	&	40.357	&	62.804	&	115.435	&	125.861	&	143.085	\\
	&	$\widehat\sigma$	&	311.332	&	273.909	&	148.766	&	120.105	&	96.564	\\\hline
Outlier 	&	$\widehat\mu$	&	119.462	&	120.844	&	122.361	&	125.893	&	143.085	\\
Deleted data	&	$\widehat\sigma$	&	129.532	&	127.406	&	125.128	&	120.009	&	96.564	\\
\hline
\end{tabular}
\label{TAB:7MDPDE_telephone_fault}
\end{table}

\begin{figure}[h]
\centering
\subfloat[$H_0$, ~ $\sigma=132$]{
\includegraphics[width=0.35\textwidth, height=0.3\textwidth] {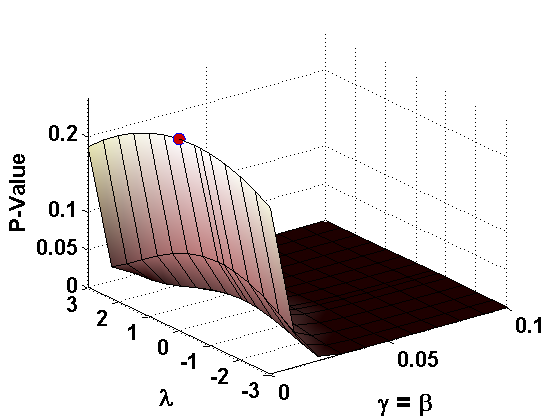}
\label{FIG:7telephoneFault_s132_test_m0}}
~ 
\subfloat[$H_0'$, ~ $\sigma=132$]{
\includegraphics[width=0.35\textwidth,  height=0.3\textwidth] {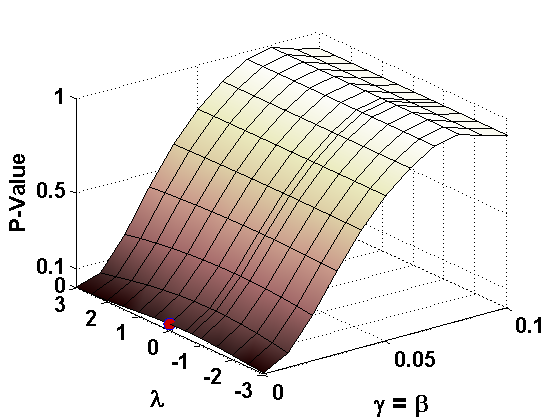}
\label{FIG:7telephoneFault_s132_test_m115}}
\\
\subfloat[$H_0$, ~ $\sigma$ unknown]{
\includegraphics[width=0.35\textwidth, height=0.3\textwidth] {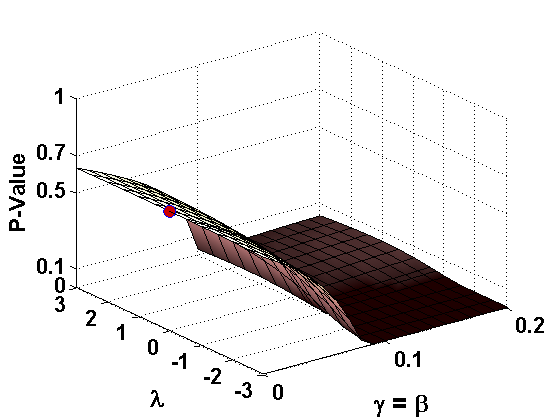}
\label{FIG:7telephoneFault_sunknown_test_m0}}
~ 
\subfloat[$H_0'$, ~ $\sigma$ unknown]{
\includegraphics[width=0.35\textwidth,  height=0.3\textwidth] {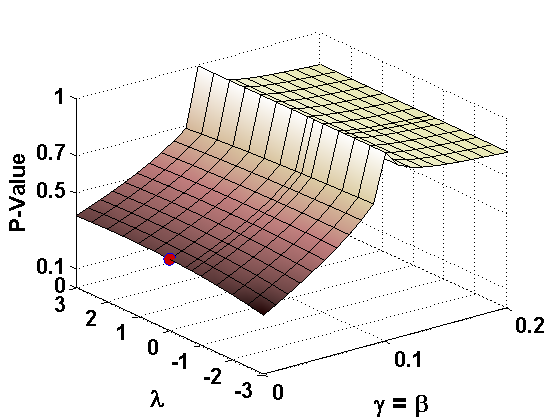}
\label{FIG:7telephoneFault_sunknown_test_m115}}
\caption{P-values for the SDT on the Telephone-fault data. The red point denotes the p-value for corresponding LRT}
 \label{FIG:7Pvalue_telephone_fault}
\end{figure}

For the present data set we consider the problem of testing two different hypothesis on the mean parameter $\mu$,
namely $H_0 :\mu = 0$ and $H_0' : \mu = 115$ against their respective  omnibus alternative. 
We consider the two cases of known and unknown $\sigma$; for the known $\sigma$ case 
we use its robust estimator $132$ for the value of $\sigma$ as suggested by \cite{Basu/etc:2013a,Basu/etc:2013b}.
However, due to the non-robust nature of the MLE,
the likelihood ratio test (and equivalently the traditional $z$-test or $t$-test) fails to reject the 
null $H_0$ but rejects the null $H_0'$ due to the presence of the large outlier; on the other hand 
after the removal of the large outlier the same test fails to reject the null $H_0'$ but soundly rejects $H_0$.
Here we apply the proposed $S$-divergence based test for these two testing problems using the 
formulas given in Ghosh et al.~(2015) and in the present paper for known and unknown $\sigma$ respectively. 
The p-values obtained by the $S$-divergence based tests (SDT) for 
the full data are presented in Figure \ref{FIG:7Pvalue_telephone_fault} along with that corresponding to LRT. 
The robust nature of the proposed SDT is clear for both the cases.



\begin{thebibliography}{}

\bibitem[\protect\citeauthoryear{Basu, Harris, Hjort, and Jones}{Basu et~al.}{1998}]{Basu/etc:1998}
Basu, A., I.~R. Harris, N.~L. Hjort, and M.~C. Jones (1998).
\newblock Robust and efficient estimation by minimising a density power
  divergence.
\newblock {\em Biometrika\/}~{\em 85}, 549--559.


\bibitem[\protect\citeauthoryear{Basu, Shioya and Park}{Basu et~al.}{2011}]{Basu/etc:2011}
Basu, A., Shioya, H. and Park, C. (2011).
\newblock {\em Statistical Inference: The Minimum Distance Approach}.
\newblock Chapman \& Hall/CRC. 

\bibitem[\protect\citeauthoryear{Basu and Lindsay}{Basu and  Lindsay}{1994}]{Basu/Lindsay:1994}
Basu, A. and B.~G. Lindsay (1994).
\newblock Minimum disparity estimation for continuous models: Efficiency,
  distributions and robustness.
\newblock {\em Annals of the Institute of Statistical Mathematics\/}~{\em 46},
  683--705.

\bibitem[\protect\citeauthoryear{Basu, Mandal, Martin, and Pardo}{Basu  et~al.}{2013a}]{Basu/etc:2013a}
Basu, A., Mandal, A., Martin, N., and Pardo, L.~(2013a).
\newblock Testing statistical hypotheses based on the density power divergence.
\newblock {\em Annals of the Institute of Statistical Mathematics\/}~{\em 65},
  319--348.
  
\bibitem[\protect\citeauthoryear{Basu, Mandal, Martin, and Pardo}{Basu et~al.}{2013b}]{Basu/etc:2013b}
Basu, A., Mandal, A., Martin, N., and Pardo, L.~(2013b).
\newblock Density Power Divergence Tests for Composite Null Hypotheses. {\em Pre-print}.


\bibitem[\protect\citeauthoryear{Cressie and Read}{Cressie and Read}{1984}]{Cressie/Read:1984}
Cressie, N. and T. R. C. Read (1984). 
\newblock Multinomial goodness-of-fit tests. 
\newblock  {\em Journal of Royal Statistical Society - B\/}~{\em 46}, 440--464.


%
%
%

\bibitem[\protect\citeauthoryear{Ghosh}{Ghosh}{2014a}]{Ghosh:2014a}
Ghosh, A.~(2014a). 
\newblock Influence Function of the Restricted Minimum Divergence Estimators : A General Form.
\newblock ArXiv Pre-print.

\bibitem[\protect\citeauthoryear{Ghosh}{Ghosh}{2014b}]{Ghosh:2014b}
Ghosh, A.~(2014b). 
\newblock Asymptotic Properties of Minimum $S$-Divergence Estimator for Discrete Models.
\newblock {\em Sankhya A: Indian Journal of Statistics\/}.


\bibitem[\protect\citeauthoryear{Ghosh and Basu}{Ghosh and Basu}{2014}]{Ghosh/Basu:2014}
Ghosh, A., Basu A. (2014). 
\newblock The Minimum S-Divergence Estimator in Continuous Models: The Basu-Lindsay Approach.
\newblock ArXiv Pre-print.


\bibitem[\protect\citeauthoryear{Ghosh, Basu and Pardo}{Ghosh et al.}{2015}]{Ghosh/etc:2015}
Ghosh, A., Basu, A., Pardo, L. (2015). 
\newblock On the Robustness of a Divergence based Test of Simple Statistical Hypotheses.
\newblock {\em J. Stat. Plann. Inf.\/}


\bibitem[\protect\citeauthoryear{Ghosh, Harris, Maji, Basu and Pardo}{Ghosh et~al.}{2013a}]{Ghosh/etc:2013a}
Ghosh, A., Harris, I. R., Maji, A., Basu, A., Pardo, L. (2013). 
\newblock A Generalized Divergence for Statistical Inference. 
\newblock {\it Technical Report}, {\bf BIRU/2013/3}, 
Indian Statistical Institute, Kolkata, India.


\bibitem[\protect\citeauthoryear{Ghosh, Maji, and Basu}{Ghosh et~al.}{2013b}]{Ghosh/etc:2013}
Ghosh, A., A.~Maji, and A.~Basu (2013).
\newblock {\em Robust Inference Based on Divergences in Reliability Systems}.
\newblock In {\em Applied Reliability Engineering and Risk Analysis. Probabilistic Models and Statistical Inference},  
Ilia Frenkel, Alex, Karagrigoriou, Anatoly Lisnianski \& Andre Kleyner, Eds, 
Dedicated to the Centennial of the birth of Boris Gnedenko, Wiley, New York, USA. 



\bibitem[\protect\citeauthoryear{Hampel, Ronchetti, Rousseeuw, and Stahel}{Hampel et~al.}{1986}]{Hampel/etc:1986}
Hampel, F.~R., E.~Ronchetti, P.~J. Rousseeuw, and W.~Stahel (1986).
\newblock {\em Robust Statistics: The Approach Based on Influence Functions}.
\newblock Wiley, New York.




\bibitem[\protect\citeauthoryear{Heritier and Ronchetti}{Heritier and Ronchetti}{1994}]{Heritier/Ronchetti:1994}
Heritier, S. and Ronchetti, E. (1994). 
\newblock Robust bounded-influence tests in general parametric models. 
\newblock  {\em J. American Stat. Ass.\/}~{\em 89}, 897--904.

\bibitem[\protect\citeauthoryear{Huber-Carol}{Huber-Carol}{1970}]{Huber/Carol:1970}
Huber-Carol, C. (1970).
\newblock {\em Etude asymptotique de tests robustes}.
\newblock Ph.\ D. thesis, ETH, Zurich.




%

\bibitem[\protect\citeauthoryear{Lehmann}{Lehmann}{1983}]{Lehmann:1983}
Lehmann, E.~L. (1983).
\newblock {\em Theory of Point Estimation}.
\newblock John Wiley \& Sons.
  
  



\bibitem[\protect\citeauthoryear{Neyman and Pearson}{Neyman and  Pearson}{1928}]{Neyman/Pearson:1928}
Neyman, J. and E.~S. Pearson (1928).
\newblock On the use and interpretation of certain test criteria for purposes
  of statistical inference.
\newblock {\em Biometrika\/}~{\em 20A}.


%
%
%
%
%

\bibitem[\protect\citeauthoryear{Rousseeuw and Ronchetti}{Rousseeuw and Ronchetti}{1979}]{Rousseeuw/Ronchetti:1979}
Rousseeuw,~P.~J. and Ronchetti,~E. (1979).
\newblock The influence curve for tests. 
\newblock {\em Research Report} {\bf 21}, Fachgruppe f\"{u}r Statistik, ETH, Zurich.


\bibitem[\protect\citeauthoryear{Rousseeuw and Ronchetti}{Rousseeuw and Ronchetti}{1981}]{Rousseeuw/Ronchetti:1981}
Rousseeuw,~P.~J. and Ronchetti,~E. (1981).
\newblock Influence curves for general statistics. 
\newblock {\em Journal of Computational and Applied Mathematics\/}~{\em 7}, 161--166.



\bibitem[\protect\citeauthoryear{Simpson}{Simpson}{1989}]{Simpson:1989b}
Simpson, D.~G. (1989).
\newblock {H}ellinger deviance test: efficiency, breakdown points, and  examples.
\newblock {\em J. American Stat. Ass.\/}~{\em 84},  107--113.


\bibitem[\protect\citeauthoryear{Toma and Broniatowski}{Toma and Broniatowski}{2011}]{Toma/Broniatowski:2010}
Toma, A. and M.~Broniatowski (2011).
\newblock Dual divergence estimators and tests: robustness results.
\newblock {\em Journal of Multivariate Analysis\/}~{\em 102}, 20--36.


\bibitem[\protect\citeauthoryear{Welch}{Welch}{1987}]{Welch:1987}
Welch, W.~J. (1987).
\newblock Rerandomizing the median in matched-pair designs.
\newblock {\em Biometrika\/}~{\em 74}, 609--614.


\bibitem[\protect\citeauthoryear{Wilks}{Wilks}{1938}]{Wilks:1938}
Wilks, S.~S. (1938).
\newblock The large sample distribution of the likelihood ratio for testing  composite hypothesis.
\newblock {\em Annals of Mathematical Statistics\/}~{\em 9}, 60--62.

\end{thebibliography}
\end{document}